%% file: pgdfreeenergy.tex
\documentclass[11pt]{article}

\usepackage{amsthm,amsmath}
\usepackage{graphicx,url,latexmac,hyperref,enumerate}
\usepackage[latin1]{inputenc}
\usepackage{ae,amsfonts,amssymb}
\theoremstyle{plain}
\newtheorem{theo}{Theorem}
\newtheorem{coro}[theo]{Corollary}

\theoremstyle{remark}

\theoremstyle{definition}

\title{Directed polymer in random environment and last passage percolation\thanks{the author acknowledges the support of the French Ministry of Education through the ANR BLAN07-2184264 grant}}

\author{Philippe Carmona\thanks{
 Laboratoire Jean Leray, UMR 6629
 Universit{\'e} de Nantes, BP 92208,
 F-44322 Nantes Cedex 03}\\
{\small \url{http://www.math.sciences.univ-nantes.fr/~carmona}}}

\newcommand{\QQ}{\mathbf{Q}}
\newcommand{\as}{a.s.}

\newcommand{\qprob}[1]{\mathbf{Q}\etp{#1}}

\renewcommand{\PP}{\mathbf{P}}
\newcommand{\Vlam}{V^{(\lambda)}}
\newcommand{\vlam}{v^{(\lambda)}}
\newcommand{\ilam}{I^{(\lambda)}}
\newcommand{\hnoe}{H_n(\omega,\eta)}
\newcommand{\fnhnoe}{\frac{H_n}{n}(\omega,\eta)}
\begin{document}
\maketitle

\begin{abstract}
The sequence of random probability measures $\nu_n$ that gives a path of length $n$, $\unsur{n}$ times the sum of the random weights collected along the paths, is shown to satisfy a large deviations principle with good rate function the Legendre transform of the free energy of the associated directed polymer in a random environment.

Consequences on the asymptotics of the typical  number of paths whose collected weight is above a fixed proportion are then drawn.
\end{abstract}

\bigskip


\textbf{Mathematic Classification} : 60K37

\newpage

\input{introduction}

\input{pfthm}

\bibliographystyle{amsplain}
\bibliography{polymeres,poly}
\end{document}

%% file: introduction.tex
\section{Introduction}
\subsection*{Last passage percolation} To each site $(k,x)$ of $\N\times \Z^d$ is assigned a random weight $\eta(k,x)$. The $(\eta(k,x))_{k\ge 1,x\in\Z^d}$ are taken IID under the probability measure $\QQ$.

The set of oriented paths of length $n$ starting from the origin is 
$$ \Omega_n =\ens{\omega=(\omega_0,\ldots,\omega_n) : \omega_i \in \Z^d, \omega_0=0, \valabs{\omega_i -\omega_{i-1}} =1}\,.$$
The weight (energy, reward) of a path is the sum of weights of visited sites:
$$ H_n = H_n(\omega,\eta) = \sum_{k=1}^n \eta(k,\omega_k)\quad(n\ge 1, \omega \in \Omega_n).$$

Observe that when $\eta(k,x)$ are Bernoulli($p$) distributed
$$ \qprob{\eta(k,x)=1} = 1 -\qprob{\eta(k,x)=0} = p \in (0,1)\,,$$
the quantity $\frac{H_n}{n}(\omega,\eta)$ is the proportion of \emph{open} sites visited by $\omega$, and it is natural to consider for $0< \rho<1$,
$$N_n(\rho) = \text{number of paths of length $n$ such that $H_n(\omega,\eta)\ge n\rho$}.$$
The problem of $\rho$-percolation, as we learnt it from Comets, Popov and Vachkovskaia~\cite{CPV07} and Kesten and Sidoravicius~\cite{kestensid07}, is to study the behaviour of $N_n(\rho)$ for large $n$ and different values of $\rho$.

\subsection*{Directed polymer in a random environment}
We are going to consider fairly general environment distributions, by requring first that they have exponential moments of any order:
$$ \lambda(\beta) = \log \qprob{e^{\beta \eta(k,x)}} < + \infty \quad(\beta \in \R)\,,$$
and second that they satisfy a logarithmic Sobolev inequality (see e.g.~\cite{MR1845806}): in particular we can apply our result to bounded support and Gaussian environments.

The polymer measure is the random probability measure defined on the set of oriented paths of length $n$ by:
$$\mu_n(\omega) = (2d)^{-n} \;\frac{e^{\beta H_n(\omega,\eta)}}{Z_n(\beta)}\qquad(\omega \in\Omega_n),,$$
with $Z_n(\beta)$ the partition function
$$ Z_n(\beta)=Z_n(\beta,\eta) = (2d)^{-n} \sum_{\omega \in \Omega_n} e^{\beta H_n(\omega,\eta)} = \prob{e^{\beta \hnoe}}\,,$$
where $\PP$ is the law of simple random walk on $\Z^d$ starting from the origin.

Bolthausen~\cite{B89} proved the existence of a deteministic limiting free energy
$$ p(\beta) = \lim_{n\to +\infty} \unsur{n} \qprob{\log Z_n(\beta)} = \QQ\, a.s. \lim_{n\to +\infty} \unsur{n}\log Z_n(\beta)\,.$$

Thanks to Jensen's inequality, we have the upper bound $p(\beta) \le \lambda(\beta)$ and it is  conjectured (and partially proved, see~\cite{CY06,CH06}) that the behaviour of a typical path under the polymer measure is diffusive iff $\beta \in \Crond_\eta$ the critical region
$$ \Crond_\eta = \ens{\beta \in \R: p(\beta) = \lambda(\beta)}\,.$$
In dimension $d=1$, $\Crond_\eta =\ens{0}$ and in dimensions $d\ge 3$, $\Crond_\eta$ contains a neighborhood of the origin (see \cite{B89,MR2249671}).

\subsection*{The main theorem}
The connection between Last passage percolation and Directed polymer in random environment is made by the family $(\nu_n)_{n\in\N}$ of random probability measures on the real line:
$$\nu_n(A) = \unsur{\valabs{\Omega_n}} \sum_{\omega \in\Omega_n} \un{\fnhnoe\in A} = \prob{\fnhnoe \in A}\,.$$
Indeed,
$$ N_n(\rho) = \sum_{\omega \in\Omega_n} \un{\hnoe \ge n\rho} = (2d)^n \nu_n([\rho,+\infty))\,.$$
The main result of the paper is
\begin{theo}
  $\QQ$ almost surely, the family $(\nu_n)_{n\in\N}$ satisfies a large deviations principle with good rate function $I=p^*$ the Legendre transform of the 
 free energy of the directed polymer.
\end{theo}

Let $m=\qprob{\eta(k,x)}$ be the average weight of a path $m=\qprob{\fnhnoe}$. It is natural to consider the quantities:
\begin{equation*}
N_n(\rho) = 
\begin{cases}
 \sum_{\omega\in\Omega_n} \un{\hnoe \ge n \rho}&\text{if $\rho\ge m$}\,,\\
 \sum_{\omega\in\Omega_n} \un{\hnoe \le n \rho}&\text{if $\rho< m$}\,.\\
\end{cases}
\end{equation*}
A simple exchange of limits $\beta\to \pm \infty$, and $n\to +\infty$, yields the following
$$ \rho^\pm = \QQ\,a.s.\lim_{n\to +\infty} \max_{\omega \in\Omega_n} \pm \fnhnoe = \lim_{\beta \to +\infty} \frac{p(\pm \beta)}{\beta} \in [0,+\infty]\,.$$
Repeating the proof of Theorem 1.1 of~\cite{CPV07} gives
\begin{coro}
  For $-\rho^- < \rho < \rho^+$, we have $\QQ$ almost surely,
$$ \lim_{n\to +\infty} (N_n(\rho))^{\unsur{n}} = (2d) e^{-I(\rho)}\,.$$
\end{coro}

We can then translate our knowledge of the critical region $\Crond_\eta$, into the following remark. Let 
$$ \Vrond_\eta=\ens{\rho \in\R : I(\rho) = \lambda^*(\rho)}\,.$$
In dimension $d=1$, $\Vrond_\eta = \ens{m}$ and in dimensions $d\ge 3$, $\Vrond_\eta$ contains a neighbourhood of $m$.

This means that in dimensions $d\ge 3$, the typical large deviation of $\fnhnoe$ close to its mean is the same as the large deviation of $\unsur{n}(\eta_1 + \cdots + \eta_n)$ close to its mean, with $\eta_i$ IID. There is no influence of the path $\omega$ : this gives another justification to  the name weak-disorder region given to the critical set $\Crond_\eta$.


%% file: pfthm.tex
\section{Proof of the main theorem}

Observe that for any $\beta\in\R$ we have:
\begin{equation}\label{eq:partition}
 \int e^{\beta n x} d{\nu_n}(x) = \prob{e^{\beta H_n(\omega,\eta)}} = Z_n(\beta)\quad \QQ\,\as\,.
\end{equation}
Consequently, since $e^{u} + e^{-u} \ge e^{\valabs{u}}$, we obtain for any $\beta >0$,
$$\limsup_{n\to +\infty} \unsur{n} \log  \etp{\int e^{\beta n {\valabs{x}}} d{\nu_n}(x)} \le p(\beta) + p(-\beta) < +\infty\,,$$
and the family $(\nu_n)_{n\ge 0}$ is exponentially tight (see Dembo and Zeitouni\cite{MR1619036}, or Feng and Kurtz\cite{MR2260560}). We only need to show now that for a lower semicontinuous function $I$, and for $x\in\R$
\begin{gather}
  \lim_{\delta \to 0} \liminf_{n\to \infty} \unsur{n} \log \nu_n((x-\delta,x+\delta)) = I(x)\,, \label{eq:liminf}\\
  \lim_{\delta \to 0} \limsup_{n\to \infty} \unsur{n} \log \nu_n([x-\delta,x+\delta]) = I(x).\label{eq:limsup}
\end{gather}
From  these, we shall infer that $(\nu_n)_{n\in\N}$ follows a large deviations principle with good rate function $I$. Eventually, equation \eqref{eq:partition} and
$$ \lim_{n\to\infty} \unsur{n}\log Z_n(\beta) = p(\beta)$$
will imply, by Varadhan's lemma that  $I$ and $p$ are Legendre conjuguate:
$$ I(x)=p^*(x) = \sup_{\beta \in\R} (x\beta -p(\beta))\,.$$

\medskip

The strategy of proof finds its origin in Varadahan's seminal paper\cite{MR1989232}, and has already succesfully been applied in \cite{CH04}. Let us define for $\lambda >0, x \in \Z, a \in\R$
$$\Vlam_n(x,a;\eta) = \log \PP^x\etp{e^{-\lambda \valabs{H_n(\omega,\eta)-a}}} = \Vlam(0,a; \tau_{o,x}\circ \eta)\,,$$
with $\tau_{k,x}$ the translation operator on the environment defined by :
$$ \tau_{k,x}\circ\eta(i,y)=\eta(k+i,x+y)\,,$$
and $\PP^x$ the law of simple random walk starting from $x$.

\medskip
\emph{Step 1} The functions $\vlam_n(a)=\qprob{\Vlam(0,a;\eta)}$ satisfy the inequality
\begin{equation}
  \label{eq:ineqvlam}
  \vlam_{n+m}(a+b) \ge \vlam_n(a) + \vlam_m(b)\qquad(n,m\in\N;\, a,b\in\R)\,.
\end{equation}
\begin{proof}
  Since $\valabs{H_{n+m}-(a+b)} \le \valabs{H_n -b} + \valabs{(H_{n+m}-H_n) -a}$ we have
\begin{align*}
  \Vlam_{n+m}(x,a;\eta) &\ge \log \PP^x\etp{e^{-\lambda \valabs{H_n -b}}
 e^{-\lambda \valabs{(H_{n+m}-H_n) -a}}} \\
&=\log\PP^x\etp{e^{-\lambda \valabs{H_n -b}} e^{\Vlam_m(0,a;\tau_{n,S_n}\circ \eta)}}\\
&= \log\sum_y \PP^x\etp{e^{-\lambda \valabs{H_n -b}} \un{S_n=y}} e^{\Vlam_m(0,a;\tau_{n,y}\circ \eta)}\\
&= \Vlam_n(x,b;\eta) + \log\etp{\sum_y \sigma_n(y)  e^{\Vlam_m(0,a;\tau_{n,y}\circ \eta)}}\\
&\ge \Vlam_n(x,b;\eta) + \sum_y \sigma_n(y)  \Vlam_m(0,a;\tau_{n,y}\circ \eta)&\text{(\small Jensen's inequality)}\,,\\
\end{align*}
with $\sigma_n$ the probability measure on $\Z^d$:
$$\sigma_n(y) = \unsur{\Vlam_n(x,b;\eta)} \PP^x\etp{e^{-\lambda \valabs{H_n-b}} \un{S_n=y}}\quad(y\in\Z^d)\,.$$

Observe that the random variables $\sigma_n(y)$ are measurable with respect to the sigma field $\Grond_n=\sigma(\eta(i,x): i\le n,x\in\Z^d)$, whereas the random variables  $\Vlam_m(0,a;\tau_{n,y}\circ \eta)$ are independent from $\Grond_n$. Hence, by stationarity,
\begin{align*}
  \vlam_{n+m}(x,a;\eta) &=\qprob{\Vlam_{n+m}(x,a;\eta)}\\
&\ge \vlam_n(b) + \sum_y \qprob{\sigma_n(y)}\qprob{\Vlam_m(0,a;\tau_{n,y}\circ \eta)} \\
&=  \vlam_n(b) + \sum_y \qprob{\sigma_n(y)}\vlam_m(a)\\
&= \vlam_n(b) + \vlam_m(a)\qprob{\sum_y \sigma_n(y)}\\
&= \vlam_n(b) + \vlam_m(a)\,.
\end{align*}
\end{proof}

\emph{Step 2} There exists a function $\ilam:\R\to\R^+$ convex, non negative, Lipschitz with constant $\lambda$, such that
\begin{equation}
  \label{eq:subadd}
  -\lim_{n\to \infty}\unsur{n} \vlam_n(a_n) = \ilam(\xi)\qquad(\text{if } \frac{a_n}{n}\to \xi\in\R))\,.
\end{equation}
\begin{proof}
  This is a standard subadditivity argument (see e.g. Varadhan~\cite{MR1989232} or Alexander~\cite{MR1428498}) combined with the Lipschitz property of $\Vlam$: from $\valabs{H_n-a} \le \valabs{H_n -b} + \valabs{a-b}$ we infer that
$$\Vlam_n(0,a;\eta)  \ge \Vlam_n(0,a;\eta) + \lambda \valabs{a-b}\,.$$
\end{proof}

\medskip
\emph{Step 3} $\QQ$ almost surely, for any $\xi\in\R$, 
\begin{equation}\label{eq:pgdps}
\lim_{n\to\infty} -\unsur{n} \log \prob{e^{-\lambda \valabs{H_n-a_n}}}=\ilam(\xi)\,.
\end{equation}
\begin{proof}
  Since the functions are Lipschitz, it is enough to prove that for any fixed $\xi\in\Q$, \eqref{eq:pgdps} holds  a.s. This is where we use the restrictive assumptions made on the distribution of the environment. If the distribution of $\eta$ is with bounded support, or Gaussian, or more generally satisfies a logarithmic Sobolev inequality, then it has the gaussian concentration of measure property (see\cite{MR1845806}): for any $1$-Lipschitz function $F$ of independent  random variables distributed as $\eta$, 
$$ \prob{\valabs{F -\prob{F}}\ge r} \le 2 e^{-r^2/2} \quad(r>0).$$
It is easy to prove, as in Proposition 1.4 of \cite{MR1939654}, that the function
$$(\eta(k,x),k\le n,\valabs{x}\le n) \to \log \prob{e^{-\lambda \valabs{H_n(\omega,\eta) -a}}}$$
is Lipschitz, with respect to the euclidean norm, with Lipschitz constant at most $\lambda \sqrt{n}$. Therefore, the Gaussian concentration of measure yields
$$ \qprob{\valabs{\Vlam_n(0,a;\eta)-\vlam_n(a)}\ge u} \le 2 e^{- \frac{\lambda^2 u^2}{2 n}}\,.$$
We conclude by a Borel Cantelli argument combined with~(\ref{eq:subadd})

\end{proof}

Observe that for fixed $\xi\in\R$, the function $\lambda \to \ilam(\xi)$ is increasing ; we shall consider the limit:
$$ I(\xi) = \lim_{\lambda \uparrow +\infty} \uparrow \ilam(\xi)$$
which is by construction non negative, convex and lower semi continuous.

\medskip
\emph{Step 4} The function $I$ satisfy (\ref{eq:liminf}) and (\ref{eq:limsup}).

\begin{proof}
  Given, $\xi \in\R$ and  $\lambda >0,\delta>0$, we have
$$ \prob{\valabs{\frac{H_n}{n}(\omega, \eta) - \xi}\le \delta} = \prob{e^{-\lambda n \valabs{\frac{H_n}{n}(\omega, \eta) - \xi}} \ge e^{-\lambda n \delta}}
 \le e^{\lambda n \delta} \prob{e^{-\lambda\valabs{H_n -n \xi}}}\,.$$
Therefore,
\begin{gather*}
\limsup \unsur{n} \log \nu_n([\xi-\delta,\xi+\delta]) \le \lambda \delta - \ilam(\xi)\\
\limsup_{\delta \to 0} \limsup \unsur{n} \log \nu_n([\xi-\delta,\xi+\delta]) \le -\ilam(\xi)
\end{gather*}
and we obtain by letting $\lambda \to +\infty$,
$$ \limsup_{\delta \to 0} \limsup \unsur{n} \log \nu_n([\xi-\delta,\xi+\delta]) \le -I(\xi)\,.$$

\medskip

Given $\xi\in\R$ such that $I(\xi)<+\infty$, and $\delta>0$,  we have for $\lambda >0$,
$$ \prob{\valabs{\frac{H_n}{n} -\xi} < \delta} \ge \prob{e^{-\lambda \valabs{H_n -n \xi}}} - e^{-\lambda \delta n}\,.$$
Hence, if we choose $\lambda>0$ large enough such that $\lambda \delta > I(\xi) \ge \ilam(\xi)$, we obtain
$$ \liminf_{n\to +\infty} \unsur{n}\log \nu_n((\xi-\delta,\xi+\delta)) \ge -\ilam(\xi)\ge -I(\xi)$$
and therefore
$$ \liminf_{\delta\to 0}\, \liminf_{n\to +\infty} \unsur{n} \log \nu_n((\xi-\delta,\xi+\delta)) \ge -I(\xi)\,.$$

\end{proof}

%% file: pgdfreeenergy.bbl
\providecommand{\bysame}{\leavevmode\hbox to3em{\hrulefill}\thinspace}
\providecommand{\MR}{\relax\ifhmode\unskip\space\fi MR }
\providecommand{\MRhref}[2]{%
  \href{http://www.ams.org/mathscinet-getitem?mr=#1}{#2}
}
\providecommand{\href}[2]{#2}
\begin{thebibliography}{10}

\bibitem{MR1428498}
Kenneth~S. Alexander, \emph{Approximation of subadditive functions and
  convergence rates in limiting-shape results}, Ann. Probab. \textbf{25}
  (1997), no.~1, 30--55. \MR{MR1428498 (98f:60203)}

\bibitem{MR1845806}
C{\'e}cile An{\'e}, S{\'e}bastien Blach{\`e}re, Djalil Chafa{\"{\i}}, Pierre
  Foug{\`e}res, Ivan Gentil, Florent Malrieu, Cyril Roberto, and Gr{\'e}gory
  Scheffer, \emph{Sur les in\'egalit\'es de {S}obolev logarithmiques},
  Panoramas et Synth\`eses [Panoramas and Syntheses], vol.~10, Soci\'et\'e
  Math\'ematique de France, Paris, 2000, With a preface by Dominique Bakry and
  Michel Ledoux. \MR{MR1845806 (2002g:46132)}

\bibitem{B89}
Erwin Bolthausen, \emph{{A note on the diffusion of directed polymers in a
  random environment.}}, Commun. Math. Phys. \textbf{123} (1989), no.~4,
  529--534.

\bibitem{MR1939654}
Philippe Carmona and Yueyun Hu, \emph{On the partition function of a directed
  polymer in a {G}aussian random environment}, Probab. Theory Related Fields
  \textbf{124} (2002), no.~3, 431--457. \MR{MR1939654 (2003m:60286)}

\bibitem{CH04}
\bysame, \emph{{Fluctuation exponents and large deviations for directed
  polymers in a random environment.}}, Stochastic Processes Appl. \textbf{112}
  (2004), no.~2, 285--308.

\bibitem{CH06}
\bysame, \emph{{Strong disorder implies strong localization for directed
  polymers in a random environment.}}, ALEA \textbf{2} (2006), 217--229.

\bibitem{CY06}
F.~Comets and N.~Yoshida, \emph{{Directed polymers in random environment are
  diffusive at weak disorder.}}, Annals of Probability \textbf{34} (2006),
  no.~5, 1746--1770.

\bibitem{CPV07}
Francis Comets, Serguei Popov, and Marina Vachkovskaia, \emph{The number of
  open paths in an oriented $\rho$-percolation model}, preprint, 2007.

\bibitem{MR2249671}
Francis Comets and Vincent Vargas, \emph{Majorizing multiplicative cascades for
  directed polymers in random media}, ALEA Lat. Am. J. Probab. Math. Stat.
  \textbf{2} (2006), 267--277 (electronic). \MR{MR2249671 (2007k:60318)}

\bibitem{MR1619036}
Amir Dembo and Ofer Zeitouni, \emph{Large deviations techniques and
  applications}, second ed., Applications of Mathematics (New York), vol.~38,
  Springer-Verlag, New York, 1998. \MR{MR1619036 (99d:60030)}

\bibitem{MR2260560}
Jin Feng and Thomas~G. Kurtz, \emph{Large deviations for stochastic processes},
  Mathematical Surveys and Monographs, vol. 131, American Mathematical Society,
  Providence, RI, 2006. \MR{MR2260560}

\bibitem{kestensid07}
Harry Kesten and Vladas Sidoravivius, \emph{A problem in last-passage
  percolation}, preprint, 2007.

\bibitem{MR1989232}
S.~R.~S. Varadhan, \emph{Large deviations for random walks in a random
  environment}, Comm. Pure Appl. Math. \textbf{56} (2003), no.~8, 1222--1245,
  Dedicated to the memory of J\"urgen K. Moser. \MR{MR1989232 (2004d:60073)}

\end{thebibliography}
